\documentclass{article}
\usepackage{amsmath,amssymb}
\newtheorem{df}{Definition}[section]
\newtheorem{thm}[df]{Theorem}
\newtheorem{lem}[df]{Lemma}
\newtheorem{prop}[df]{Proposition}

\title{Spin weighted wavelets on the sphere}
\author{I. Iglewska-Nowak\footnote{West Pomeranian University of Technology in Szczecin, School of Mathematics, al. Piast\'ow 17, 70--310 Szczecin, Poland}}

\begin{document}

\maketitle

\bibliographystyle{amsplain}

\begin{abstract}In the present paper, a construction of spin weighted spherical wavelets is presented. It is based on approximate identities, the wavelets are defined for a continuous set of parameters, and the wavelet transform is invertible directly by an integral.\end{abstract}

\begin{bfseries}Keywords:\end{bfseries} spherical wavelets, spin weighted functions \\
\begin{bfseries}MSC2010:\end{bfseries} 42C40, 83C60

\section{Introduction}

In the paper~\cite{NP66} Newman and Penrose introduced an invariant differential operator~$\eth$ on the sphere and used it to define spin $s$ spherical harmonics. Spin weight of a function is defined by considering rotations around the radial direction of unit vectors tangent to the sphere, understood as vectors in the complex line bundle, tangent to the sphere at each point. An $\mathcal L^2$-function~${}_sf$ is said to have spin weight $s\in\mathbb Z$, ${}_sf\in\mathcal L_s^2$, if it transforms by ${}_sf\mapsto\exp(is\psi)\,{}_sf$, whenever a tangent vector at a point transforms by $\mathbf m\mapsto \exp(i\psi)\,\mathbf m$. In this sense, it is a tensorlike quantity~\cite{NP66}. It can be understood as a function of a point on the sphere \emph{and} a direction in the tangent space. However, a spin weighted function is given as a scalar function of the two spherical variables~$\vartheta\in[0,\pi]$ and $\varphi\in[0,2\pi)$.

We refer the reader to the original paper~\cite{NP66}, as well as to~\cite{GMNRS67}, where the relationship between spin weighted spherical harmonics  and elements of the representation matrices of the rotation group~$SO(3)$ is pointed out. An exhaustive explanation (with a slightly other machinery) and many applications can be found in the excellent textbook~\cite{TC03}, further examples of the usefulness of this notion in physics are presented in~\cite{TC07}.

In the recent years, several wavelet constructions for spin weighted functions were presented. One of them are the \emph{needlet-type spin wavelets} \cite{GM10,GM10a}, a generalization of the spherical needlets \cite{GM09a,GM09b}, further developed to the \emph{mixed needlets} \cite{GM11}. Spin wavelets are given as series of $s$-zonal harmonics~${}_sZ_l$,
$$
{}_sZ_l(x)=\sum_{m=-l}^l{}_s^{}Y_l^m(x)\,\overline{{}_s^{}Y_l^m(\hat e)},
$$
where $\hat e$ denotes the north pole of the sphere. A discretized set of wavelets constitutes a nearly tight frame and the wavelet transform is inverted by the frame methods.

The construction of the mixed needlets is similar. The main difference is the choice of the basic functions
$$
\sum_{m=-l}^l{}_sY_l^m(x)\,\overline{Y_l^m(\hat e)},
$$
'mixed' from the spin $s$ weighted and the spin $0$ weighted (i.e., scalar) spherical harmonics. Further, the scale set is a priori discrete. Again, the reconstruction results from the frame theory and is given as a series.

The \emph{directional spin wavelets}~\cite{MLBPW17} are intended to probe the directional intensity of spin signals. The theory is based on group theoretical approach and the difference between this theory and that based on approximate identities is discussed in~\cite{IIN15CWT} (for the case of spin $0$ functions). The scale coefficient is a priori discrete but there is a direct integral reconstruction formula.

There is one more construction of a wavelet transform for spin weighted spherical functions worth to be mentioned. In~\cite{SMB09} spin $2$ functions are transformed into spin $0$ functions and they are a subject  of the undecimated isotropic wavelet transform.

In the present paper we introduce spin wavelet analysis being a generalization of the spherical wavelet analysis based on approximate identities. Both the scale and position parameters are continuous and a direct reconstruction is possible by an integral.

The paper is organized as follows. In Section~\ref{sec:sphere}, basic information about spin weighted functions are collected. A definition of the spin weighted wavelets and the wavelet transform, as well as proofs of its invertibility and isometry are the content of Section~\ref{sec:bwt}. An example of spin weighted wavelets is presented in Section~\ref{sec:example}, and Section~\ref{sec:discussion} contains a discussion on the relevance of this research.

\section{Preliminaries}\label{sec:sphere}

\subsection{Spin-weighted functions on the sphere}

Let $\mathcal S$ denote the $2$-dimensional unit sphere in the Euclidean space~$\mathbb{R}^3$. Points on the sphere are given in the spherical coordinates $(\vartheta,\varphi)\in[0,\pi]\times[0,2\pi)$,
$$
\mathcal S\ni x=\left(\begin{array}{c}\sin\vartheta\cos\varphi\\\sin\vartheta\sin\varphi\\\cos\vartheta\end{array}\right).
$$
Unless it leads to misunderstandings, the points are identified with their coordinates. $d\sigma$ denotes the Lebesgue measure of~$\mathcal S$,  $d\sigma(\vartheta,\varphi)=\sin\vartheta\,d\vartheta\,d\varphi$. The scalar product of $f,g\in\mathcal L^2(\mathcal S)$ is defined by
$$
\left<f,g\right>_{\mathcal L^2(\mathcal S)}=\int_{\mathcal S}\overline{f(x)}\,g(x)\,d\sigma(x),
$$
such that $\|f\|_2^2=\left<f,f\right>$.

The spherical harmonics of degree $l\in\mathbb N_0$ and order~$m=-l,-l+1,\dots,l-1,l$ are defined by
\begin{align*}
Y_l^m(\vartheta,\varphi)&=(-1)^m\sqrt{\frac{2l+1}{4\pi}\frac{(l-m)!}{(l+m)!}}\cdot P_l^m(\cos\vartheta)\cdot \exp(im\varphi)&(m\geq0),\\
Y_l^m&=\overline{Y_l^{-m}}&(m<0),
\end{align*}
where $P_l^m$ denotes the associated Legendre polynomial. The spherical harmonics constitute an orthonormal basis for~$\mathcal L^2(\mathcal S)$ and they are functions of spin weight~$0$.

Denote by $\mathcal L_s^2(\mathcal S)$ the space of square integrable functions over the sphere with spin weight~$s$. The spin-raising operator
$$
\eth\colon\;\mathcal L_s^2(\mathcal S)\to\mathcal L_{s\!+\!1}^2(\mathcal S)
$$
and the spin-lowering operator
$$
\bar\eth\colon\;\mathcal L_{s\!+\!1}^2(\mathcal S)\to\mathcal L_s^2(\mathcal S)
$$
are given by
\begin{align*}
\eth\,{}_sf(\vartheta,\varphi)&=-\left(\partial_\vartheta+\frac{i}{\sin\vartheta}\,\partial_\varphi-s\cot\vartheta\right){}_sf(\vartheta,\varphi)\\
   &=-\sin^s\vartheta\cdot\left(\partial_\vartheta+\frac{i}{\sin\vartheta}\,\partial_\varphi\right)\,\frac{{}_sf(\vartheta,\varphi)}{\sin^{s}\vartheta},\\
\bar\eth\,{}_sf(\vartheta,\varphi)&=-\left(\partial_\vartheta-\frac{i}{\sin\vartheta}\,\partial_\varphi+s\cot\vartheta\right){}_sf(\vartheta,\varphi)\\
   &=-\frac{1}{\sin^{s}\vartheta}\cdot\left(\partial_\vartheta-\frac{i}{\sin\vartheta}\,\partial_\varphi\right)\left[{}_sf(\vartheta,\varphi)\cdot\sin^s\vartheta\right].
\end{align*}

The spin weighted spherical harmonics are defined as
\begin{align*}
_s^{}Y_l^m&=\sqrt{\frac{(l-s)!}{(l+s)!}}\,\,\eth^s\,Y_l^m,&0\leq s\leq l,\\
_s^{}Y_l^m&=\sqrt{\frac{(l+s)!}{(l-s)!}}\,\,(-1)^s\,\,\bar\eth^{-s}\,Y_l^m,&-l\leq s\leq0,
\end{align*}
and for each $s\in\mathbb Z$ they are an orthonormal basis for~$\mathcal L_s^2(\mathcal S)$.
Every $\mathcal L_s^2(\mathcal S)$-function~${}_sf$ has a unique representation as the Fourier series
$$
{}_sf=\sum_{l=|s|}^\infty\sum_{m=-l}^l\,{}_s^{}\widehat f_l^m\,{}_s^{}Y_l^m
$$
with the Fourier coefficients ${}_s^{}\widehat f_l^m=\left<{}_s^{}Y_l^m,f\right>_{\mathcal L^2}$, convergent in $\mathcal L^2$-norm. The space spanned by $\{{}_s^{}Y_l^m\colon\,m=-l,-l+1,\dots,l\}$ is denoted by~${}_s\mathcal H_l$. 

The following holds for the spin weighted spherical harmonics:
$$
\sum_{l=|s|}^\infty\sum_{m=-l}^l{}_s^{}Y_l^m(\vartheta_1,\varphi_1)\,\overline{{}_s^{}Y_l^m(\vartheta_2,\varphi_2)}
   =\delta(\cos\vartheta_1-\cos\vartheta_2)\,\delta(\varphi_1-\varphi_2),
$$
cf. \cite[formula (3.15)]{GMNRS67}. Additionally, consider the kernel
\begin{equation}\label{eq:kernel_K}
{}_s\mathcal K_l(x,y):=\sum_{m=-l}^l{}_s^{}Y_l^m(x)\cdot\overline{{}_s^{}Y_l^m(y)}.
\end{equation}
According to~\cite[formula~(2.59)]{TC03} it is equal to
\begin{equation}\label{eq:sKl_as_sY}
{}_s\mathcal K_l(x,y)=(-1)^{-s}\,
\sqrt\frac{2l+1}{4\pi}\cdot{}_s^{}Y_l^{-s}(\vartheta_3,\varphi_3)\cdot \exp(-is\chi_3)
\end{equation}
where $x=(\vartheta_1,\varphi_1)$, $y=(\vartheta_2,\varphi_2)$ in the spherical coordinates, and $(\varphi_3,\vartheta_3,\chi_3)$ are the Euler angles of the composition of rotations with the Euler angles $(0,-\vartheta_2,-\varphi_2)$ and $(\varphi_1,\vartheta_1,0)$. We shall prove the reproducing property of~${}_s\mathcal K_l$.

\begin{prop}\label{prop:reproducing_kernel} ${}_s\mathcal K_l$ is the reproducing kernel of ${}_s\mathcal H_l$, i.e.,
$$
({}_s\mathcal K_l\ast {}_sf)(x):=\int_{\mathcal S^2}{}_s\mathcal K_l(x,y)\cdot {}_sf(y)\,d\sigma(y)={}_sf(x)
$$
for each ${}_sf\in{}_s\mathcal H_l$. The operator ${}_sf\mapsto{}_s\mathcal K_l\ast{}_sf$ is the projection~${}_sf_l$ from $\mathcal L_s^2$ onto~${}_s\mathcal H_l$.
\end{prop}

\begin{bfseries}Proof. \end{bfseries} Let an ${}_s\mathcal H_{l^\prime}$-function
$$
{}_sf=\sum_{m=-l^\prime}^{l^\prime}{}_s^{}\widehat f_{l^\prime}^m\,{}_s^{}Y_{l^\prime}^m
$$
be given. Its convolution with ${}_s\mathcal K_l$ equals
\begin{align*}
({}_s\mathcal K_l&\ast {}_sf)(x)=\int_{\mathcal S^2}{}_s\mathcal K_l(x,y)\cdot{}_sf(y)\,d\sigma(y)\\
&=\sum_{m=-l}^l\sum_{m^\prime=-l^\prime}^{l^\prime}{}_s^{}\widehat f_{l^\prime}^{m^\prime}{}_s^{}Y_l^m(x)
   \int_{\mathcal S^2}\overline{{}_s^{}Y_l^m(y)}\cdot{}_s^{}Y_{l^\prime}^{m^\prime}(y)\,d\sigma(y)\\
&=\sum_{m=-l}^l\sum_{m^\prime=-l^\prime}^{l^\prime}{}_s^{}\widehat f_{l^\prime}^{m^\prime}{}_s^{}Y_l^m(x)\,\delta_{ll^\prime}\,\delta_{mm^\prime}
   ={}_sf(x)\,\delta_{ll^\prime}.
\end{align*}
For an $\mathcal L_s^2$-function note that $\mathcal L_s^2=\bigoplus_{l=|s|}^\infty{}_s\mathcal H_l$.\hfill$\Box$

Next we show that the kernel~$_s\mathcal K_l$ is uniformly bounded.

\begin{prop}\label{prop:estimation_kernel}The estimation
\begin{equation}\label{eq:sKl_estimation}
\left|{}_s\mathcal K_l(x,y)\right|\leq C\cdot l^{2|s|+1}
\end{equation}
holds uniformly in~$x$ and~$y$ with a positive constant~$C$.
\end{prop}

The following auxiliary result will be used in the proof of the Proposition.

\begin{lem}\label{lem:estimation_P0k}
For $t\in[-1,1]$, the Jacobi polynomials satisfy
$$
\left|(1+t)^k\cdot P_n^{(0,k)}(t)\right|\leq2^k.
$$
\end{lem}

\begin{bfseries}Proof. \end{bfseries}For $k=0$, the Jacobi polynomials reduce to the Legendre polynomials, which are uniformly bounded by~$1$,
$$
\left|P_n^{(0,0)}(t)\right|=\left|P_n(t)\right|\leq1\qquad\text{for }t\in[-1,1],
$$
see \cite[formulae~(8.962.2) and~(8.917.7)]{GR}. Further, \cite[formula and~(8.961.6)]{GR} yields
$$
\left|(1+t)\cdot P_n^{(0,k+1)}(t)\right|\leq\frac{n+k+1}{n+\tfrac{k}{2}+1}\cdot\left| P_n^{(0,k)}(t)\right|+\frac{n+1}{n+\tfrac{k}{2}+1}\cdot\left| P_{n+1}^{(0,k)}(t)\right|.
$$
The assertion is obtained by induction.\hfill$\Box$

\begin{bfseries}Proof of Proposition~\ref{prop:estimation_kernel}. \end{bfseries}Formula~\eqref{eq:sKl_as_sY} together with
$$
{}_s^{}Y_l^m(\vartheta,\varphi)=(-1)^m\cdot\sqrt\frac{2l+1}{4\pi}\cdot d_{-m,s}^j(\vartheta)\cdot \exp(im\varphi),
$$
(compare~\cite[second formula on p.~54]{TC03}) yields
\begin{equation}\label{eq:kernel_as_dlss}
\left|{}_s\mathcal K_l(x,y)\right|=\frac{2l+1}{4\pi}\cdot\left|d_{s,s}^l(\vartheta)\right|,
\end{equation}
$d_{m,k}^l(\vartheta)$ being the elements of the Wigner d-matrices. According to~\cite[formula~(99)]{gA13},
\begin{equation}\label{eq:dlmk_as_Legenedre}
d_{m,k}^l(\vartheta)=\frac{1}{2^k}\sqrt{\frac{(l-k)!}{(l-m)!}\frac{(l+k)!}{(l+m)!}}\cdot\sin^{k-m}\vartheta\cdot(1+\cos\vartheta)^m\cdot P_{l-k}^{(k-m,k+m)}(\cos\vartheta),
\end{equation}
and the symmetry relation~\cite[formula~(93)]{gA13},
\begin{equation}\label{eq:dlmk_symmetry}
d_{m,k}^l(\vartheta)=(-1)^{m-k}\cdot d_{-m,-k}^l(\vartheta),
\end{equation}
$d_{s,s}^l$ can be written as
$$
d_{s,s}^l(\vartheta)=\frac{(1+\cos\vartheta)^{|s|}}{2^{|s|}}\cdot P_{l-|s|}^{(0,2|s|)}(\cos\vartheta)
$$
where $P_n^{(\alpha,\beta)}$ is the Jacobi polynomial. Thus, by Lemma~\ref{lem:estimation_P0k},
\begin{equation}\label{eq:estimation_dlss_near0}
\left|d_{s,s}^l(\vartheta)\right|\leq\left(\frac{2}{1+\cos\vartheta}\right)^{|s|}.
\end{equation}
for $\vartheta\in[0,\pi)$.  On the other hand, by~\cite[formula~(94)]{gA13} for integer indices $l$, $m$, $k$,
$$
d_{m,-k}^l(\vartheta)=(-1)^{l-m}\cdot d_{m,k}^l(\pi-\vartheta),
$$
together with~\eqref{eq:dlmk_symmetry} and~\eqref{eq:dlmk_as_Legenedre},
\begin{align*}
d_{s,s}^l(\pi-\vartheta)&=(-1)^{l-s}\,d_{-|s|,|s|}^l(\vartheta)
   =\frac{(-1)^{l-s}}{2^{|s|}}\cdot\frac{\sin^{2|s|}\vartheta}{(1+\cos\vartheta)^{|s|}}\cdot P_{l-|s|}^{(2|s|,0)}(\cos\vartheta)\\
&=(-1)^{l-s}\cdot\sin^{2|s|}\frac{\vartheta}{2}\cdot P_{l-|s|}^{(2|s|,0)}(\cos\vartheta).
\end{align*}
In order to estimate this expression we use \cite[Main Theorem, p.~980]{FW85} with $m=1$. Coefficient $A_0$ is equal to~$1$, see \cite[p.~994]{FW85}.
$$
\sin^{2|s|}\frac{\vartheta}{2}\,P_n^{(2|s|,0)}(\cos\vartheta)=\frac{\Gamma(n+2|s|+1)}{\Gamma(n+1)}\,\sqrt\frac{\vartheta}{\sin\vartheta}
   \,\left[\frac{J_{2|s|}(N\vartheta)}{N^{2|s|}}+\mathcal O(N^{-1})\right],
$$
where
$$
N=n+|s|+\frac{1}{2},\qquad\vartheta\in(0,\pi),
$$
and~$J_{2|s|}$ denotes the Bessel function. The $\mathcal O$-term is uniform with respect to $\vartheta\in[0,\pi-\epsilon]$. Since
\begin{align*}
\frac{\Gamma(n+2|s|+1)}{\Gamma(n+1)}&=\Pi_{\iota=n+1}^{n+2|s|}\,\iota=\Pi_{\iota=1}^{|s|}(n+\iota)(n+2|s|+1-\iota)\\
&\leq\Pi_{\iota=1}^{|s|}(n+|s|+\tfrac{1}{2})^2=(n+|s|+\tfrac{1}{2})^{2|s|},
\end{align*}
and by \cite[formula~(43)]{lL98},
$$
|J_{2|s|}(t)|<\frac{b}{\sqrt[3]{2|s|}}
$$
uniformly in~$t$ with $b=0.674885\dots$, the absolute value of $d_{s,s}^l(\pi-\vartheta)$ is for $\vartheta\in(0,\pi-\epsilon]$ bounded by
$$
\left|d_{s,s}^l(\pi-\vartheta)\right|\leq(l+\tfrac{1}{2})^{2|s|}\cdot\sqrt\frac{\vartheta}{\sin\vartheta}\cdot\left(\frac{1}{\sqrt[3]{2|s|}}+\frac{C}{l+\tfrac{1}{2}}\right),
$$
where~$C$ is a constant depending on~$\epsilon$. The limiting case $\vartheta\to0$ is obtained by continuity and convergence of $\vartheta/\sin\vartheta$. Thus,
$$
\left|d_{s,s}^l(\pi-\vartheta)\right|\leq C\cdot(l+\tfrac{1}{2})^{2|s|}
$$
uniformly for $\vartheta\in[0,\pi-\epsilon]$ and with another constant~$C$. This, together with~\eqref{eq:estimation_dlss_near0}, yields
$$
\left|d_{s,s}^l(\vartheta)\right|\leq C\cdot l^{2|s|}
$$
uniformly in~$\vartheta\in[0,\pi]$ for a constant~$C>0$. The assertion follows by~\eqref{eq:kernel_as_dlss}.\hfill$\Box$

\subsection{Rotations of spin weighted functions and their zonal product}\label{subs:rotation}

$SO(3)$-rotations~$\Upsilon$ are given in the Euler angles, $\Upsilon=(\alpha,\beta,\gamma)$, denoting the rotation of the points on the sphere in a fixed coordinate system about $OZ$, $OY$, and $OZ$ axes by $\gamma$, $\beta$, and $\alpha$, respectively. The rotation of a spin weighted function~${}_sf$ on the sphere is defined by
\begin{equation}\label{eq:rotation_function}
(\mathcal R_\Upsilon\,{}_sf)(x):=\exp(-is\kappa)\cdot{}_sf\left(\Upsilon^{-1}x\right),\qquad \Upsilon\in SO(3),
\end{equation}
with $\kappa\in[0,2\pi)$ being the third Euler angle of $\Upsilon^{-1}\Xi$, where $\Xi$ is the rotation given by the Euler angles $(\varphi,\vartheta,0)$ for $x=(\vartheta,\varphi)$ in the spherical coordinates. The exponential factor $\exp(-is\kappa)$ ensures that the rotation of~${}_sf$ is again a function with spin weight~$s$ \cite[Subsection~2.4]{MLBPW17} \cite[Appendix~A]{CL05}.

According to \cite{MLBPW17,CL05}, the spin weighted spherical harmonics satisfy
\begin{equation}\label{eq:rotation_Yslm}
\mathcal R_\Upsilon\,{}_s^{}Y_l^m=\sum_{k=-l}^lD_l^{km}(\Upsilon)\cdot{}_s^{}Y_l^k,
\end{equation}
where $D_l^{km}$, $l\in\mathbb N_0$, $k,m\in\mathbb Z$, $|k|,|m|\leq l$, are the Wigner D-functions, i.e., the matrix elements of the irreducible unitary representation of the rotation group $SO(3)$, orthogonal to each other in the sense that
\begin{equation}\label{eq:scalar_product_D}\begin{split}
\left<D_{l_1}^{k_1m_1},D_{l_2}^{k_2m_2}\right>_{\mathcal L^2\left(SO(3)\right)}&:=\int_{SO(3)}\overline{D_{l_1}^{k_1m_1}(\Upsilon)}\,D_{l_2}^{k_2m_2}(\Upsilon)\,d\nu(\Upsilon)\\
&=\frac{8\pi^2}{2l_1+1}\,\delta_{l_1l_2}\delta_{k_1k_2}\delta_{m_1m_2},
\end{split}\end{equation}
$d\nu$ being the Haar measure on $SO(3)=\mathcal S\times\mathcal T$, $\mathcal T=[0,2\pi)$, normalized such that $\int_{SO(3)}d\nu(\Upsilon)=8\pi^2$.

Analogously to the spin $0$case \cite{EBCK09}, we introduce the zonal product of spin weighted square integrable functions by
$$
({}_{s_1}\!f\,\hat\ast\,{}_{s_2}g)(x,y)=\int_{SO(3)}(\mathcal R_\Upsilon\,{}_{s_1}f)(x)\cdot(\mathcal R_\Upsilon\,{}_{s_2}g)(y)\,d\nu(\Upsilon),\qquad x,y\in\mathcal S.
$$
It is well-defined by the H\"older inequality and
$$
\int_{SO(3)}\left|(\mathcal R_\Upsilon\,{}_sf)(x)\right|^2\,d\nu(\Upsilon)=2\pi\int_{\mathcal S}|{}_sf(x)|^2\,d\sigma(x)<\infty.
$$

\begin{lem}\label{lem:zp_as_series} For ${}_sf,\,{}_sg\in\mathcal L_s^2$,
$$
({}_sf\,\hat\ast\,\overline{{}_sg})(x,y)
   =\sum_{l=|s|}^\infty\frac{8\pi^2}{2l+1}\sum_{m=-l}^{l}{}_s^{}\widehat f_l^m\cdot\overline{{}_s^{}\widehat g_l^m}\cdot{}_s\mathcal K_l(x,y)\qquad\text{a.e.}
$$

\end{lem}

\begin{bfseries}Proof. \end{bfseries}For an integer $N>|s|$, set
\begin{align*}
\mathfrak r_N(x,y)=\int_{SO(3)}&\sum_{l_1=|s|}^N\sum_{m_1=-l_1}^{l_1}\sum_{l_2=|s|}^N\sum_{m_2=-l_2}^{l_2}
   {}_s^{}\widehat f_{l_1}^{\,m_1}\cdot\overline{{}_s^{}\widehat g_{l_2}^{m_2}}\\
&\cdot(\mathcal R_\Upsilon\,{}_s^{}Y_{l_1}^{m_1})(x)\cdot\overline{(\mathcal R_\Upsilon\,{}_s^{}Y_{l_2}^{m_2})(y)}\,d\nu(\Upsilon).
\end{align*}
By the H\"older inequality and the  square summability of the Fourier coefficients ${}_s\widehat f_l^m$, ${}_s\widehat g_l^m$,
\begin{equation}\label{eq:rN_limit}
\mathfrak r_N\longrightarrow{}_sf\,\hat\ast\,{}_sg\qquad\text{for}\qquad N\to\infty
\end{equation}
in $\mathcal L^2$-sense.
Equation~\eqref{eq:rotation_Yslm} yields
\begin{align*}
\mathfrak r_N&(x,y)
=\sum_{l_1\!=\!|s|}^N\sum_{m_1\!=\!-\!l_1}^{l_1}\sum_{l_2\!=\!|s|}^N\sum_{m_2\!=\!-\!l_2}^{l_2}
   {}_s^{}\widehat f_{l_1}^{m_1}\cdot\overline{{}_s^{}\widehat g_{l_2}^{m_2}}\\
&\cdot\sum_{k_1\!=\!-\!l_1}^{l_1}\sum_{k_2\!=\!-\!l_2}^{l_2}{}_s^{}Y_{l_1}^{k_1}(x)\cdot\overline{{}_s^{}Y_{l_2}^{k_2}(y)}
   \cdot\int_{SO(3)}D_{l_1}^{k_1m_1}(\Upsilon)\cdot \overline{D_{l_2}^{k_2m_2}(\Upsilon)}\,d\nu(\Upsilon).
\end{align*}
Thus, by~\eqref{eq:scalar_product_D},
\begin{align*}
\mathfrak r_N(x,y)&=\sum_{l_1\!=\!|s|}^N\sum_{m_1\!=\!-\!l_1}^{l_1}\sum_{l_2\!=\!|s|}^N\sum_{m_2\!=\!-\!l_2}^{l_2}
   {}_s^{}f_{l_1}^{m_1}\cdot\overline{g_{l_2}^{m_2}}\\
&\cdot\sum_{k_1\!=\!-\!l_1}^{l_1}\sum_{k_2\!=\!-\!l_2}^{l_2}{}_s^{}Y_{l_1}^{k_1}(x)\cdot\overline{{}_s^{}Y_{l_2}^{k_2}(y)}
   \cdot\frac{8\pi^2}{2l_1+1}\,\delta_{l_1l_2}\delta_{k_1k_2}\delta_{m_1m_2}\\
&=\sum_{l=|s|}^N\frac{8\pi^2}{2l+1}\sum_{m=-l}^{l}{}_s^{}f_l^m\cdot\overline{g_l^m}\cdot\sum_{k=-l}^l{}_sY_l^k(x)\cdot\overline{{}_sY_l^k(y)}.
\end{align*}
Using notation~\eqref{eq:kernel_K}, we can write it as
$$
\mathfrak r_N(x,y)=\sum_{l=|s|}^N\frac{8\pi^2}{2l+1}\sum_{m=-l}^{l}{}_s^{}f_l^m\cdot\overline{g_l^m}\cdot{}_s\mathcal K_l(x,y).
$$
The assertion follows by~\eqref{eq:rN_limit}.\hfill$\Box$

\section{The continuous wavelet transform of spin weighted functions}\label{sec:bwt}

The following definitions generalize the concept of the bilinear wavelets and the bilinear wavelet transforms based on approximate identities \cite{FGS-book,EBCK09,IIN15CWT} to the case of spin weighted functions. 

\begin{df}\label{def:bilinear_wavelets} Let $\alpha:\mathbb R_+\to\mathbb R_+$ be a weight function, piecewise continuous and bounded on compact sets. A family $\{{}_s\Psi_\rho\}_{\rho\in\mathbb R_+}\subseteq\mathcal L_s^2(\mathcal S)$ is called an admissible spherical wavelet of spin weight~$s$ if it satisfies the following conditions:
\begin{enumerate}
\item for $l\in\mathbb{N}_0$, $l\geq|s|$,
\begin{equation}\label{eq:admbwv1}
\sum_{m=-l}^{l}\int_0^\infty\left|{}_s^{}(\widehat{\Psi_\rho})_l^m\right|^2\cdot\alpha(\rho)\,d\rho=\frac{2l+1}{8\pi^2},
\end{equation}
\item for a fixed~$\rho$,
\begin{equation}\label{eq:admbwv2}
\sum_{l=|s|}^\infty l^{2|s|}\cdot\sum_{m=-l}^{l}\left|{}_s^{}(\widehat{\Psi_\rho})_l^m\right|^2<\infty.
\end{equation}
\end{enumerate}
\end{df}

\begin{df}\label{def:bilinear_wt} Let $\{{}_s\Psi_\rho\}_{\rho\in\mathbb R_+}$ be a spherical wavelet with spin weight~$s$. Then, the spherical wavelet transform
$$
\mathcal W_\Psi\colon\mathcal L_s^2(\mathcal S)\to\mathcal L^2(\mathbb R_+\times SO(3))
$$
is defined by
\begin{equation}\label{eq:bwt}
(\mathcal W_\Psi{}_s f)(\rho,\Upsilon)=\int_{\mathcal S}\overline{(\mathcal R_\Upsilon\,{}_s\Psi_\rho)(x)}\cdot{}_sf(x)\,d\sigma(x).
\end{equation}
\end{df}

The integral in~\eqref{eq:bwt} is absolutely convergent, since both the wavelet and the analyzed function are square integrable.

\begin{thm}\label{thm:inversion}If $\{{}_s\Psi_\rho\}_{\rho\in\mathbb R_+}$ is an admissible wavelet, then the wavelet transform is invertible by
\begin{equation}\label{eq:bwt_synthesis}
{}_sf(x)=\lim_{R\to 0}\int_R^{1/R}\!\!\int_{SO(3)}(\mathcal R_\Upsilon\,{}_s\Psi_\rho)(x)\cdot(\mathcal W_\Psi\,{}_sf)(\rho,\Upsilon)\,d\nu(\Upsilon)\cdot\alpha(\rho)\,d\rho
\end{equation}
with limit in $\mathcal L^2$-sense.
\end{thm}

\begin{bfseries}Proof. \end{bfseries}Denote the inner integral in~\eqref{eq:bwt_synthesis} by $\mathfrak r_\rho(x)$. By~\eqref{eq:bwt}, it is equal to
$$
\mathfrak r_\rho(x)
   =\int_{SO(3)}\int_{\mathcal S}(\mathcal R_\Upsilon\,{}_s\Psi_\rho)(x)\cdot\overline{(\mathcal R_\Upsilon\,{}_s\Psi_\rho)(y)}\cdot{}_sf(y)\,d\sigma(y)\,d\nu(\Upsilon).
$$
Since the functions ${}_s\Psi_\rho$ and~${}_sf$ are square integrable, the integral is convergent and the order of integration may be changed according to Fubini's theorem,
\begin{align*}
\mathfrak r_\rho(x)&=\int_{\mathcal S}\int_{SO(3)}
   (\mathcal R_\Upsilon\,{}_s\Psi_\rho)(x)\cdot\overline{(\mathcal R_\Upsilon\,{}_s\Psi_\rho)(y)}\,d\nu(\Upsilon)\cdot {}_sf(y)\,d\sigma(y)\\
&=\int_{\mathcal S}\left({}_s\Psi_\rho\,\hat\ast\,\overline{{}_s\Psi_\rho}\right)(x,y)\cdot {}_sf(y)\,d\sigma(y).
\end{align*}
According to Lemma~\ref{lem:zp_as_series}, the zonal product of the wavelets can be expressed as
\begin{equation}\label{eq:zonal_product_PsiOmega}
\left({}_s\Psi_\rho\,\hat\ast\,\overline{{}_s\Psi_\rho}\right)(x,y)=\sum_{l=|s|}^\infty\frac{8\pi^2}{2l+1}
   \sum_{m=-l}^{l}\left|{}_s^{}(\widehat{\Psi_\rho})_l^m\right|^2\cdot{}_s\mathcal K_l(x,y).
\end{equation}
By~\eqref{eq:sKl_estimation} and~\eqref{eq:admbwv2}, the series on the right-hand-side of~\eqref{eq:zonal_product_PsiOmega} is uniformly bounded. Since ${}_sf\in\mathcal L^2(\mathcal S)\subset\mathcal L(\mathcal S)$, the order of integration over~$\mathcal S$ and summation with respect to~$l$ in
$$
\mathfrak r_\rho(x)=\int_\mathcal S\,\sum_{l=|s|}^\infty\frac{8\pi^2}{2l+1}
   \sum_{m=-l}^{l}\left|{}_s^{}(\widehat{\Psi_\rho})_l^m\right|^2\cdot{}_s\mathcal K_l(x,y)\cdot{}_sf(y)\,d\sigma(y)
$$
can be changed. With Proposition~\ref{prop:reproducing_kernel} we obtain
\begin{align*}
\mathfrak r_\rho&=\sum_{l=|s|}^\infty\frac{8\pi^2}{2l+1}\sum_{m=-l}^{l}\left|{}_s^{}(\widehat{\Psi_\rho})_l^m\right|^2\cdot{}_sf_l\\
&=\sum_{l=|s|}^\infty\sum_{k=-l}^l\left(\frac{8\pi^2}{2l+1}\sum_{m=-l}^{l}\left|{}_s^{}(\widehat{\Psi_\rho})_l^m\right|^2\right)\cdot{}{}_s^{}\widehat f_l^k\cdot{}_s^{}Y_l^k.
\end{align*}
The coefficients
$$
{}_s^{}(\widehat{\mathfrak r_\rho})_l^k:=\left(\frac{8\pi^2}{2l+1}\sum_{m=-l}^{l}\left|{}_s^{}(\widehat{\Psi_\rho})_l^m\right|^2\right)\cdot{}{}_s^{}\widehat f_l^k
$$
are summable by~\eqref{eq:admbwv2} and boundedness of~${}_s^{}\widehat f_l^k$. Thus, they are square summable and $\mathfrak r\in\mathcal L_s^2(\mathcal S)$. This justifies the usage of notation~${}_s^{}(\widehat{\mathfrak r_\rho})_l^k$. For any~$R$, the $\mathcal L^2$-norm $L(R)$ of
$$
f-\int_R^{1/R}\mathfrak r_\rho\cdot\alpha(\rho)\,d\rho
$$
satisfies
\begin{equation}\label{eq:L(R)}
[L(R)]^2=\sum_{l=|s|}^\infty\sum_{k=-l}^l
   \left[1-\int_R^{1/R}\frac{8\pi^2}{2l+1}\sum_{m=-l}^{l}\left|{}_s^{}(\widehat{\Psi_\rho})_l^m\right|^2\cdot\alpha(\rho)\,d\rho\right]^2\cdot\left({}_s^{}\widehat f_l^k\right)^2.
\end{equation}
By~\eqref{eq:admbwv1}, for each $(l,k)$ the difference in brackets on the right-hand-side of~\eqref{eq:L(R)} is bounded. Thus, by square summability of~$\left({}_s^{}\widehat f_l^k\right)_{l,k}$, the series in~\eqref{eq:L(R)} is uniformly bounded. Consequently,
\begin{align*}
\lim_{R\to0}&[L(R)]^2\\
&=\sum_{l=|s|}^\infty\sum_{k=-l}^l\lim_{R\to0}
   \left[1-\int_R^{1/R}\frac{8\pi^2}{2l+1}\sum_{m=-l}^{l}\left|{}_s^{}(\widehat{\Psi_\rho})_l^m\right|^2\cdot\alpha(\rho)\,d\rho\right]^2\cdot\left({}_s^{}\widehat f_l^k\right)^2\\
&=0,
\end{align*}
and the assertion follows.\hfill$\Box$\\[-1em]

Similarly as in the case of spin $0$ functions~\cite{IIN15CWT}, the wavelet transform is an isometry.

\begin{thm}\label{thm:isometry} Let $\{{}_s\Psi_\rho\}_{\rho\in\mathbb R_+}$ be an admissible spherical wavelet of spin weight~$s$ and ${}_sf,{}_sg\in\mathcal L_s^2(\mathcal S)$. Then,
$$
\left<\mathcal W_\Psi\,{}_sf,\mathcal W_\Psi\,{}_sg\right>=\left<{}_sf,{}_sg\right>
$$
where the scalar product in the wavelet phase space is given by
$$
\left<F,G\right>_{\mathcal L^2(\mathbb R_+\times SO(3))}=\int_0^\infty\!\!\int_{SO(3)}\overline{F(\rho,\Upsilon)}\,G(\rho,\Upsilon)\,d\nu(\Upsilon)\,\alpha(\rho)\,d\rho.
$$
\end{thm}

The proof is analogous to that of Theorem~\ref{thm:inversion}.

\begin{bfseries}Proof.\end{bfseries} For a fixed~$\rho$, Denote by~$\mathfrak r_\rho$ the integral
\begin{align*}
\mathfrak r_\rho&=\int_{SO(3)}\overline{(\mathcal W_\Psi\,{}_sf)(\rho,\Upsilon)}\cdot(\mathcal W_\Psi\,{}_sg)(\rho,\Upsilon)\,d\nu(\Upsilon)\\
&=\int_{SO(3)}\int_{\mathcal S}(\mathcal R_\Upsilon\,{}_s\Psi_\rho)(x)\,\overline{{}_sf(x)}\,d\sigma(x)
   \cdot\int_{\mathcal S}\overline{(\mathcal R_\Upsilon\,{}_s\Psi_\rho)(y)}\,{}_sg(y)\cdot d\sigma(y)\,d\nu(\Upsilon).
\end{align*}
A change of the integration order yields
\begin{align*}
\mathfrak r_\rho=\int_{\mathcal S}\int_{\mathcal S}&
   \underbrace{\int_{SO(3)}(\mathcal R_\Upsilon\,{}_s\Psi_\rho)(x)\cdot\overline{(\mathcal R_\Upsilon\,{}_s\Psi_\rho)(y)}\cdot d\nu(\Upsilon)}
   _{{}_s\Psi_\rho\,\hat\ast\,\overline{{}_s\Psi_\rho}(x,y)}\\
&\cdot\alpha(\rho)\,d\rho\cdot\overline{{}_sf(x)}\,d\sigma(x)\cdot{}_sg(y)\,d\sigma(y).
\end{align*}
According to~\eqref{eq:zonal_product_PsiOmega}, this expression is equal to
$$
\mathfrak r_\rho=\int_{\mathcal S}\int_{\mathcal S}\sum_{l=|s|}^\infty\frac{8\pi^2}{2l+1}\sum_{m=-l}^l\left|{}_s^{}(\widehat{\Psi_\rho})_l^m\right|^2\cdot{}_s\mathcal K_l(x,y)
   \cdot\overline{{}_sf(x)}\,d\sigma(x)\cdot{}_sg(y)\,d\sigma(y).
$$
By~\eqref{eq:sKl_estimation} and~\eqref{eq:admbwv2}, the order of integration and summation may be changed,
\begin{align*}
\mathfrak r_\rho&=\int_{\mathcal S}\sum_{l=|s|}^\infty\frac{8\pi^2}{2l+1}\sum_{m=-l}^l\left|{}_s^{}(\widehat{\Psi_\rho})_l^m\right|^2
   \cdot\int_{\mathcal S}{}_s\mathcal K_l(x,y)\cdot\overline{{}_sf(x)}\,d\sigma(x)\cdot{}_sg(y)\,d\sigma(y)\\
&=\int_{\mathcal S}\sum_{l=|s|}^\infty\frac{8\pi^2}{2l+1}\sum_{m=-l}^l\left|{}_s^{}(\widehat{\Psi_\rho})_l^m\right|^2\cdot\overline{{}_sf_l(y)}\cdot{}_sg(y)\,d\sigma(y).
\end{align*}
Since
$$
{}_sf_l(y)\leq\|{}_s\mathcal K_l(\circ,y)\|_{\mathcal L^2}\cdot\|{}_sf\|_{\mathcal L^2}\leq C\cdot l^{2|s|+1},
$$
the series is absolutely convergent, and the order of summation and integration may be changed,
$$
\mathfrak r_\rho=\sum_{l=|s|}^\infty\frac{8\pi^2}{2l+1}\sum_{m=-l}^l\left|{}_s^{}(\widehat{\Psi_\rho})_l^m\right|^2
   \cdot\int_{\mathcal S}\overline{{}_sf_l(y)}\cdot{}_sg(y)\,d\sigma(y).
$$
By the orthogonality of the spin weighted spherical harmonics,
$$
\int_{\mathcal S}\overline{{}_sf_l(y)}\cdot{}_sg(y)\,d\sigma(y)=\sum_{k=-l}^l\overline{{}_s^{}\widehat f_l^k}\cdot{}_s^{}\widehat g_l^k.
$$
The assertion follows by the square summability of the Fourier coefficients ${}_s^{}\widehat f_l^k$, ${}_s^{}\widehat g_l^k$ and condition~\eqref{eq:admbwv1}, with the same arguments as in the proof of Theorem~\ref{thm:inversion}.\hfill $\Box$\\

\begin{bfseries}Remark.\end{bfseries} In the case of spin $0$ functions, there is a distinguished class of zonal functions, i.e., those invariant with respect to the rotation about the $z$-axis. For zonal wavelets, $SO(3)$-rotations simplify to $\mathcal S$-translations and the wavelet coefficients depend on the scale and the position (instead of rotation)~\cite{IIN15CWT,IIN15PW}. We do not attempt to generalize the notion of zonal wavelets to the spin weighted case. There exists a family of spin weighted functions invariant with respect to the rotations about the $z$-axis, namely those with Fourier coefficients ${}_s\widehat f_l^m$ that do not vanish only for $m=0$. However, a pure coordinate rotation of a spin weighted function is not a function with the same spin order, see~\cite[footnote~7]{MLBPW17}, and in the rotation formula~\eqref{eq:rotation_function} the exponential factor $\exp(-is\kappa)$ must occur. Thus, the resulting formulae would be as sophisticated as those defining the wavelet transform for the whole rotation group.

\section{Example of spin weighted wavelets}\label{sec:example}

It is straightforward to construct spin weighted wavelets with a generating function, analogous to the wavelets of spin $0$ such as Abel-Poisson, Gauss-Weierstrass, or Poisson wavelets. A generalization of these examples is presented in \cite[Theorem~2.6]{IIN16FDW}, and the following theorem gives their adaptation to the spin weighted case.

\begin{thm}A family $\{{}_s\Psi_\rho\}$ of $\mathcal L_s^2$-functions given by
\begin{equation*}\begin{split}
{}_s^{}(\widehat \Psi_\rho)_l^0&=\frac{1}{2\pi}\cdot\sqrt\frac{2l+1}{2}\cdot\sqrt\frac{4^c\cdot a}{\Gamma(2c)}
   \cdot\left(\rho^a\,[q_\gamma(l)]^b\right)^c\,\exp\left(-\rho^a\,[q_\gamma(l)]^b\right),\\
{}_s^{}(\widehat \Psi_\rho)_l^k&=0,\qquad k\ne0,
\end{split}\end{equation*}
for $l\in\mathbb N_0$, 
where~$q_\gamma$ is a polynomial of degree~$\gamma$, strictly positive for a positive~$l$, and $a$, $b$, $c$ --- some positive constants, is a wavelet with respect to the weight \mbox{$\alpha(\rho)=1/\rho$}.
\end{thm}

\begin{bfseries}Proof. \end{bfseries}Conditions~\eqref{eq:admbwv1} and~\eqref{eq:admbwv2} can be verified by a direct calculation.\hfill$\Box$\\

\section{Discussion}\label{sec:discussion}

It has been shown in \cite{IIN15CWT} that there exist two essentially different continuous wavelet transforms for square integrable functions over the sphere, namely that based on group theory and that based on approximate identities. The first one has been recently generalized to the case of spin weighted functions~\cite{MLBPW17}. The present paper constitutes the 'approximate identities' counterpart to the work of McEwen \emph{et al.} Significant advantages of the construction introduced here are the facts that the wavelet transform is defined for a continuous set of parameters and that it is invertible by an integral. Similarly as in~\cite{MLBPW17}, admissible wavelets can reveal directional features such that directional intensity of spin signal can be analyzed. Further, the constraints on a function family to be a wavelet are quite weak such that a wise class of functions can play a role of a wavelet. One can construct wavelets with desired featured, such as angular selectivity or uncertainty product (compare \cite{IIN17AS,IIN17USW} for the spin $0$ case).

Moreover, it seems plausible that wavelet frames exist, both continuous and discrete, compare \cite{IIN16FDW} for the scalar wavelet transform. This will be the object of my future research.\\

\end{document}